\documentclass[twocolumn]{revtex4}

\usepackage{graphicx}

\def\a{\alpha}
\def\b{\beta}

\def\la{\lambda}

\def\pa{\partial}

\def\xb{{\bf x}}

\def\o+{\oplus}

\def\ss{\subset}

\def\<{\langle}
\def\>{\rangle}

\def\({\left(}
\def\){\right)}
\def\[{\left[}
\def\]{\right]}
\def\=#1{\bar #1}

\def\.#1{\dot #1}
\def\^#1{\widehat #1}

\def\"#1{\ddot #1}

\def\eeq{\end{equation}}
\def\beq{\begin{equation}}

\def\beql#1{\begin{equation} \label{#1}}

\def\eqref#1{(\ref{#1})}

\def\qb{{\bf q}}

\def\xb{{\bf x}}

\begin{document}

\title{Higher order normal modes}

\author{Giuseppe Gaeta\footnote{giuseppe.gaeta@unimi.it}}
\affiliation{
Dipartimento di Matematica, Universit\`a degli Studi
di Milano, via Saldini 50, 20133 Milano (Italy); \\ {\rm and} \\ SMRI,  00058 Santa Marinella (Italy) }

\author{Sebastian Walcher\footnote{walcher@matha.rwth-aachen.de}} 
\affiliation{Mathematik A, RWTH Aachen, 52056 Aachen (Germany)} 

%\date{{\giorno}}

\begin{abstract}
\noindent
Normal modes are intimately related to the quadratic approximation of a potential at its hyperbolic equilibria. Here we extend the notion to the case where the Taylor expansion for the potential at a critical point starts with higher order terms, and show  that such an extension shares some of the properties of standard normal modes. Some symmetric examples are considered in detail. \\ { } \\ 
\quad \quad \quad {\it Dedicated to James Montaldi on his 25+something anniversary}
\end{abstract}

\maketitle

\section{Introduction}

The concept of \emph{normal modes} is a fundamental one in the
study of Hamiltonian dynamical systems \cite{LLM,Gant,M2,M3}; it is
based on the quadratic part of the Taylor expansion of the
Hamiltonian around a non-degenerate (isolated and hyperbolic)
stable equilibrium, and under certain fairly general assumptions
it can be conveniently employed also in considering higher order
expansions of the Hamiltonian. That is, under such assumptions (including a non-resonance condition) normal modes persist, at
least locally, when one considers also higher order terms -- or
for nonlinear Hamiltonian dynamics \cite{NLNM1,NLNM2} -- as also
studied by James in a series of papers \cite{Mon1,Mon2}.

Normal modes span the dynamics of the quadratic Hamiltonian, i.e.
any motion for this can be described as the superposition of
normal modes; more precisely, normal modes define invariant lines
in a neighborhood $U \simeq R^n$ of the equilibrium in the
position space, and these provide a basis for $U$.

Here we want to discuss the (partial) generalization of this
concept to the case where the equilibrium is stable and isolated,
but \emph{not} hyperbolic. This resonates with some recent studies
in a different field, i.e. \emph{liquid crystals}. More precisely, we will find some connection with recent studies of liquid crystals described by
higher order tensor order parameter \cite{Vir,GV1,QiV,GV2}.

We will show that normal modes -- in the sense of invariant lines
-- also exist for the dynamics of fully nonlinear homogeneous
Hamiltonian systems (more precisely, we will confine ourselves to
systems with a natural Hamiltonian $H = T + V$ and a potential
$V(x)$ homogeneous of degree $k > 2$). On the other hand, unless some very special conditions (e.g. based on the symmetry of the full system) are met, it is \emph{not} possible to describe a generic motion as a (nonlinear) superposition of these higher order normal modes.

\section{Normal modes}
\label{sec:NM}

We will consider natural Hamiltonians \beq H \ = \ T \ + \ V \ = \
\frac{{\bf p}^2}{2 \, m} \ + \ V ({\bf q}) \eeq with ${\bf q} \in
{\bf R}^n$; the potential $V$ is such that the origin is an
isolated stable equilibrium point.

Our discussion could be extended to encompass a general
Hamiltonian on a symplectic manifold $M$ of dimension $2n$, and a
non-degenerate equilibrium point $p_0$ for this; however our
considerations will be \emph{local}, so by Darboux' theorem we can
always consider a standard symplectic form.

The physical interest of the natural Hamiltonian case surely
justifies considering this special case. {Moreover this will help
keeping the discussion and notation simpler}, focusing on the key
points.\footnote{In particular, all the discussion in this
section, i.e. referring to the quadratic case, can easily be set
in terms of a general Hamiltonian $H({\bf p},{\bf q})$.}

We can then consider the Taylor expansion of $V$ around ${\bf q} =
0$; we thus have, in view of the non-degeneracy assumption,
$$ V({\qb}) \ = \ V_0 \ + \ \frac12 \, \( \qb , A \qb \) \ + \ \mathtt{h.o.t.} \ , $$
where $V_0 = V ({\bf 0})$ is a constant -- which can be set to
zero with no loss of generality -- and $A$ is a symmetric tensor
of order two, i.e. a symmetric matrix, \beql{eq:Amat} A_{ij} \ = \
\( \frac{\pa^2 V}{\pa q^i \pa q^j} \)_0 \ . \eeq

We will now truncate the potential at order two, i.e. omit the higher order terms and deal just with $V = 1/2 (\qb , A \qb)$; the relevant case for normal modes is that where the origin is a stable fixed point, hence $A$ is a positive definite matrix.

This matrix $A$ will have eigenvalues $\la_i$ and corresponding (say
normalized) eigenvectors $\phi_{(i)}$; these in turn correspond to normal
modes in a very well known way, which we briefly recall to fix
notation \cite{LLM,Gant,M2,M3}.

\begin{enumerate}

\item The eigenvectors $\phi_{(i)}$ identify \emph{invariant
lines}: if a motion has initial conditions $\qb (0) = \a (0)
\phi_{(i)}$, $\dot{\qb} (0) = \b (0) \phi_{(i)}$, then $\qb (t) =
\a (t) \phi_{(i)}$ (and hence $\dot{\qb} (t) = \b (t) \phi_{(i)}$)
for all $t$.

\item The eigenvectors $\phi_{(i)}$ ($i = 1,...,n$) span the
linear space ${\bf R}^n$;

\item A general motion, i.e. one with generic initial conditions
$\qb (0)$, $\dot{\qb} (0)$ can then be decomposed in terms of
normal modes. In fact, by the property just mentioned, there will
exist constants $\a_i (0)$, $\b_i (0)$ such that
$$ \qb (0) \ = \ \sum_{i=1}^n \a_i (0) \, \phi_{(i)} \ ; \ \ \
\dot{\qb (0)} \ = \ \sum_{i = 1}^n \b_i (0) \phi_{(i)} \ . $$ The
evolution of the system will then just be described by
$$ \qb (t) \ = \ \sum_{i=1}^n \a_i (t) \, \phi_{(i)} \ ; \ \ \
\dot{\qb (t)} \ = \ \sum_{i = 1}^n \b_i (t) \phi_{(i)} \ , $$ with
$\a_i (t)$ and $\b_i (t)$ {being} exactly as in the normal mode solution
with initial data $\{\a_i (0) , \b_i (0)\}$. In other words, any
solution for this quadratic Hamiltonian will be a \emph{linear
superposition of normal modes solutions}.

\item As a consequence of the previous item, in the case there exist several normal modes associated {with} the same frequency, the whole linear space spanned by the associated eigenvectors is made of eigenvectors and hence of normal modes with the given frequencies.

\end{enumerate}

As anticipated, we will find that properties (1) and (2) extend to
the higher order case, while properties (3) and (4), related to linearity,
do not.

\section{Eigenvectors of tensors}
\label{sec:eigenv}

{We now consider} the case where $\qb = 0$ is an isolated
critical point, but {we will drop} the hypothesis it is hyperbolic.
Actually, we want to assume that the matrix $A$ defined in
\eqref{eq:Amat} is \emph{identically zero}, i.e. that the
potential $V$ is \emph{fully nonlinear} at the considered critical
point.

Thus we have to go further in the Taylor expansion -- even if we
want to stop at the first significant term -- and actually if we
want the critical point to be a stable one it is needed to
consider a fourth order\footnote{This is quite different from the
case of liquid crystals mentioned above: there the relevant tensor
is of order three \cite{Vir,GV1,QiV,GV2}.} tensor $T$, \beql{eq:Bten}
T_{ijk\ell} \ := \ \frac{1}{24} \( \frac{\pa^4 H}{\pa q^i \pa q^j
\pa q^k \pa q^\ell} \)_{0} \ . \eeq

Thus, disregarding higher order terms, we want to consider a
potential expressed in local coordinates which is homogeneous of
order four and which has an isolated equilibrium at the origin.
{After getting} rid of inessential (additive and multiplicative)
constants, this can be written as \beq V \ = \ T_{ijk\ell} \ q^i
\, q^j \, q^k \, q^\ell \ . \eeq

The ``direct extension'' of properties (1) and (2) recalled above
for usual normal modes would require to consider \emph{eigenvectors for
the fourth order tensor $T$} (rather than the matrix $A$) and see
if these span the whole space.

The notion of eigenvector of tensors is not so well known in general, but it
can be defined and it has been studied (with a revival of interest
in recent times) both from the algebraic point of view \cite{T1,T2,T3}
and in connection with dynamics \cite{Roh1,Roh2,Roh3,Roh4,RoW,Walg,KLPW}  (see also \cite{OrbSpa,Wal}); as already mentioned, it has also been recently considered in connection with critical points of a constrained potential with applications in the Physics of Liquid Crystals \cite{Vir,GV1,QiV,GV2}.

We will not introduce and discuss the notion of eigenvectors of tensors right
away, but {we will start} with a general discussion of one dimensional eigenspaces of
homogeneous polynomial maps.

\subsection{Eigenspaces of homogeneous polynomial vector fields}

We will simultaneously discuss the real and the complex case in this section. We {first recall} a classical result from Algebraic Geometry; see \cite{Shafa} for more about this.

\medskip\noindent
{\bf Theorem 1 (Bezout).} {\it Let $\{ f_1 , ... f_m \}$ be homogeneous polynomials of degree $n$, $f_i : {\bf C}^{m} \to {\bf C}$. Then the number of common zeros of the $f_i$ in projective space ${\bf P}^{m}$ is either infinite or equal to $n^m$, counting multiplicities.}
\bigskip

For real polynomials this theorem provides only limited information about the nature of these critical points; in particular, we cannot infer how many of these are \emph{real}.\footnote{This is of course the same situation met in the Fundamental Theorem of Algebra: we know that a polynomial of degree $n$ in one variable always has $n$ roots (counting multiplicities) but, with no further study, we do not know how many of these are real.} But, since the complex conjugate of every solution is also a solution (with the same multiplicity, as can be shown), we conclude that a real solution exists whenever $n$ is odd (and there are only finitely many solutions).

We now {denote by ${\bf K}$ the} real or complex numbers and consider a homogeneous polynomial map
\begin{equation}\label{hommap}
 B:\,{\bf K}^q\to {\bf K}^q, \quad x\mapsto\left(\begin{array}{c} B^1(x)\\ \vdots\\ B^q(x)\end{array}\right)
\end{equation}
with each $B^i$ homogeneous of degree $p\geq 2$. In coordinates we have
\begin{equation}
B^i = B^i_{j_1 ... j_q} x^{j_1} ... x^{j_q}.
\end{equation}
We stipulate that the coefficients $B^i_{j_1 ... j_q}$ are symmetric with respect to permutations of $j_1,\ldots,j_q$; this choice makes the coefficients unique. It is possible to identify $B$ with its coefficients, and consider it as an element of the coefficient space (which is just some ${\bf K}^N$.) A formal definition follows next.

\medskip\noindent
{\bf Definition.} {\it Let $B$ be given as in \eqref{hommap}. Then a nonzero $v\in {\bf K}^q$ is called an {\it eigenvector} of $B$ if there exists an $\alpha\in {\bf K}$ such that $B(v)=\alpha  v$.}

\bigskip\noindent
Some remarks are in order here:

\begin{enumerate}

\item Every nonzero scalar multiple of an eigenvector $v$ is also an eigenvector. Therefore it makes sense to call ${\bf K}\,v$ an eigenspace of $B$.

\item On the other hand, the notion of eigenvalue is problematic for homogeneous maps of degree $>1$, since
$$
B(\beta v)=\beta^p\alpha v=(\beta^{p-1}\alpha)\cdot(\beta v)
$$
for any $\beta$, and thus one may replace the ``eigenvalue'' $\alpha$ by $\beta^{p-1}\alpha$.

\item One may use this property to scale ``eigenvalues'' to be either $0$ or $1$ in the complex case, and also in the real case when the degree $p$ is even; for the real case with odd $p$ one may achieve ``eigenvalue'' $0$, $1$ or $-1$ by scaling.

\item Alternatively, one may prescribe that the (Euclidean) norm of an eigenvector should be equal to $1$; then the notion of eigenvalue becomes relevant. This will be done later for tensors and their gradient systems.

\end{enumerate}

We now list a number of results on eigenspaces of homogeneous polynomial maps; most of these are based on the work of H.~Rohrl \cite{Roh1,Roh2,Roh3,Roh4,RoW}.

\medskip\noindent
{\bf Theorem 2.} {\it Let $B$ be as in \eqref{hommap}, and ${\bf K}={\bf C}$. Then the following hold.
\begin{itemize}

\item[(i)] The number of one-dimensional eigenspaces of $B$ is either infinite or equal to
$$ N_R \ = \ \frac{p^q \ - \ 1}{p -1}, $$
counting multiplicities.

\item[(ii)] Test for multiplicity one: Let $v\in{\bf C}^q$ be nonzero and $B(v)=\alpha v$ with some $\alpha\not=0$. Then ${\bf C}v$ corresponds to a solution of multiplicity one if and only if the Jacobian $DB(v)$ does not admit the eigenvalue $\alpha$.

\item[(iii)] If the equation $B(x)=0$ has only the trivial solution $x = 0$ then the number of one-dimensional eigenspaces of $B$ is finite.

\item[(iv)] There is an open and dense subset of coefficient space such that every $B$ with coefficients in this subset admits a basis for ${\bf C}^q$ of eigenvectors.

\item[(v)] There is an open and dense subset of coefficient space such that every $B$ with coefficients in this subset admits exactly $N_R$ different one-dimensional eigenspaces.

\item[(vi)] If the coefficients $B^i_{j_1 ... j_q}$ are algebraically independent over the rational number field then the equation $B(x)=0$ has only the trivial solution, and $B$ admits exactly $N_R$ different one-dimensional eigenspaces.

\end{itemize} }

\medskip\noindent
{\bf Proof.}
We just sketch some arguments for the proofs, and give references. (See also the review \cite{Wal}.) The first assertion is due to Rohrl \cite{Roh1,Roh3}, the second is derived from the familiar criterion for multiplicity one (i.e. invertibility of the linearization). The third assertion is shown e.g. in \cite{Walg}, based on the fact that any projective variety of positive dimension intersects every hyperplane. The fourth and fifth assertion go essentially back to Rohrl \cite{Roh4}, although the full statement given in this paper is not correct, and the proof has to be modified. See the Appendix of \cite{KLPW} for a full discussion. The last statement is again due to Rohrl \cite{Roh1}; the algebraic independence condition guarantees that the multiplicity one criterion is always satisfied. \hfill $\odot$
\bigskip

We turn to the real setting. For proofs and references concerning the following statements we refer to \cite{Wal}. (In some of the proofs analytic techniques enter the picture.)

\medskip\noindent
{\bf Theorem 3.}  {\it Let $B$ be as in \eqref{hommap}, and ${\bf K}={\bf R}$. Then the following hold.

\begin{itemize}

\item[(i)] If the dimension $q$ is odd then there exists a one-dimensional real eigenspace of $B$.

\item[(ii)] If the dimension $q$ is even and the degree $p$ of $B$ is even then  there exists a one-dimensional real eigenspace of $B$.

\item[(iii)] If the complexification admits finitely many one-dimensional eigenspaces, then the number of real eigenspaces is congruent to $N_R$ modulo $2$.

\end{itemize}
}

\subsection{Radial solutions of fully nonlinear dynamical systems}
\label{sec:Rohrl}

Rohrl was  interested in (one-dimensional) eigenspaces of homogeneous polynomial maps because they give rise to special solutions of an associated differential equation, similar to the linear case. Rohrl considered first order differential equations, and we paraphrase his result here (see the original work in \cite{Roh1} and \cite{Roh2}).

\medskip\noindent
{\bf Theorem 4.} {\it Let $B$ be as in \eqref{hommap}, and consider the ordinary differential equation
$$\dot x=B(x)$$
in ${\bf K}^n$. Then every one-dimensional eigenspace of $B$ is an invariant set for this differential equation.

For nonzero $v$ with $B(v)=\alpha v$, some $\alpha\in{\bf K}$, one obtains solutions with the ansatz $x(t)=\xi(t)\cdot v$, which leads to the one-dimensional equation $\dot \xi=\alpha\xi^p$.}
\bigskip

{There is a straightforward extension of this approach to second order equations.}

\medskip\noindent
{\bf Theorem 5.} {\it Let $B$ be as in \eqref{hommap}, and consider the  second order ordinary differential equation in ${\bf K}^n$
$$ \ddot{x} \ = \ B(x) \ . $$

\begin{itemize}

\item[(i)] Then every  nonzero $v$ with $B(v)=\alpha v$, some $\alpha\in{\bf K}$,  gives rise to special solutions of the differential equation, via the ansatz $x(t)=\gamma(t)\cdot v$, which leads to the one-dimensional second order equation $\ddot \gamma=\alpha\gamma^p$.

\item[(ii)] {In case ${\bf K}={\bf R}$ a phase plane analysis of the associated system $$\dot y^1=y^2,\,\dot y^2=\alpha\left(y^1\right)^p$$ yields the first integral $\psi=2\alpha \left(y^1\right)^{p+1}-(p+1)\left(y^2\right)^2$. When $\alpha\not=0$ then the level sets of $\psi$ are bounded if and only if $\alpha<0$ and $p$ is odd. In this case, the level sets are orbits of periodic solutions of the second order system; in all other cases every nonconstant solution obtained by the ansatz is unbounded.}
\end{itemize}
}

\medskip\noindent{\bf Proof} (Sketch). {Note that $0$ is the only stationary point of the system. It is elementary to see that the level sets are bounded, hence compact, if and only if $\alpha<0$ and $p$ is odd. By standard Poincar\'e-Bendixson theory for planar systems, the only possible limit sets of points on a compact level set containing more than one point are closed orbits, thus they must coincide with the level sets by connectedness. In every other case, each level set that contains more than one point is unbounded, and any $\alpha$ or $\omega$ limit point of a solution starting on such a level set consists of a stationary point, by Poincar\'e-Bendixson. Unboundedness of the solution follows. \hfill $\odot$}

\medskip\noindent
{\bf Remark 1.} It is worth noting that in the case $p=2$ the special solutions {from} the theorem correspond to a well-known class of special functions. Indeed, the second order equation
$$ \ddot{z} \ = \ \alpha \ z^2 $$
becomes, upon employing the first integral,
$$ \dot{z}^2 \ = \ \frac23 \ \alpha \ z^3 \ + \ c $$
with some constant $c$, and this is the differential equation for a Weierstrass $\wp$-function. So, elliptic functions appear in a natural manner.  \hfill $\odot$

\subsection{Critical points on the unit sphere}
\label{sec:unit}

In this section we are only interested in the real case ${\bf K}={\bf R}$. We consider a symmetric tensor $T_{i_1...i_n}$ of order $n$ on ${\bf R}^m$, and we associate with this a polynomial
$$ P_n (x) \ := \ T_{i_1 ... i_n} x^{i_1} ... x^{i_n} \ ;$$ note that here the dimension $m$  of the ambient space and the degree $n$ of the polynomial are not related. In the following, $P_n$ will also be called the \emph{potential}; it will also be just denoted as $P$, when we do not need to emphasize its degree.
Conversely, as is well known, the algebra of homogeneous polynomials of degree $n$ in ${\bf R}^m$ is isomorphic to the algebra of symmetric tensors of the same order $n$ over ${\bf R}^m$.

Consider now the gradient of $P_n$, i.e. the $m$-dimensional vector
$$ \nabla P_n \ = \ \( \frac{\pa P_n}{\pa x^1} \ , \ ... \ , \ \frac{\pa P_n}{\pa x^m} \) \ ; $$ here of course each component is a homogeneous function of degree $(n-1)$ in the $x^i$.

{We define} an {\it eigenvector of the tensor} $T$ to be an eigenvector of  $(\nabla P_n)$.

For eigenvectors of tensors we obtain an improvement of earlier results concerning the real case. We collect them in the following.

\medskip\noindent
{\bf Theorem {6}.}
{\it Let $v\in{\bf R}^m$ be of Euclidean norm one. Then

\begin{itemize}

\item[(i)] $v$ is an eigenvector of $\nabla P_n$ if and only if $v$ is a critical point of $P_n$ on the unit sphere  $S^{m-1} \ss {\bf R}^m$.

\item[(ii)] The real homogeneous gradient map $\nabla P_n$ admits a real eigenvector.

\item[(iii)] If we have finitely many critical points $x_k$, then the sum of the indices of all critical points is equal to the (Euler-Poincar\'e) characteristic $\chi (S^{m-1})$ of the ambient sphere, thus:
\beql{eq:index} \sum_k \iota (x_k) \ = \ \chi (S^{m-1}) \ = \ \cases{2 & if $m$ is odd \cr 0 & if $m$ is even \cr} \eeq

\end{itemize}}

\medskip\noindent
{\bf Proof.} To prove the first assertion $(i)$, introduce a Lagrange multiplier $\la$ and consider the modified potential
\beq \^P (\xb) \ := \ P (\xb) \ - \ \frac12 \, \la \ |\xb|^2 \ . \eeq
The gradient of $\^P$ is given by
$$ \nabla \^P \ = \ \nabla P \ - \ \la \, \xb \ ; $$
hence the solutions to $\nabla \^P = 0$ are exactly the points on the unit sphere such that $(\nabla P) (\xb_0)$ is collinear to $\xb_0$; these identify (unit length) eigenvectors of $T$ and hence eigenspaces. The second assertion  $(ii)$ is then clear since $P_n$ attains maximum and minimum on the compact unit sphere.
{The notion of \emph{index} of a critical point $a$ is defined in \cite{Milnor} via the Brouwer degree. When $a$ is nondegenerate (i.e.,  the derivative at $a$ is invertible) it is equal to the sign of its Jacobian determinant.
For the proof of the last assertion $(iii)$ see \cite{Milnor,ArnODE}.} \hfill $\odot$
\bigskip

In view of the results above we will adopt the convention that whenever reference is made to eigenvalues associated {with} eigenvectors, it is understood that these are associated {with} eigenvectors of unit length.\footnote{Note that for $n$ odd this still leaves an ambiguity, as the eigenvalues for $v$ and $-v$ differ by sign; it would actually be convenient to consider these as two distinct eigenvalues, corresponding to distinct eigenvectors, in view of the discussion in Sect. \ref{sec:unit}, see \cite{GV1,GV2}.}

We also note that (like for matrices) if $T$ depends on parameters then the eigenvectors and eigenvalues will in general depend on these parameters. We anticipate that also the \emph{number} of independent eigenvectors (that is, eigenspaces) can vary depending on such parameters. This situation is met already in the simplest nontrivial case, i.e. for completely symmetric cubic tensors in three-dimensional space\footnote{For completely symmetric cubic tensors in two-dimensional space we get a degenerate situation; see the discussion in \cite{Vir}. We will see later on, in Section \ref{sec:ExaI}, that a similar degeneration is met for quartic tensors in two-dimensional space.}. For a full discussion of this case we refer to \cite{GV1} (with a more physical approach) and especially {to} \cite{GV2}.

The situation is specially simple when the ambient space is ${\bf R}^3$ (a case of clear physical interest!) and hence we work in $S^2 \ss {\bf R}^3$ and all critical points are non-degenerate. {In this case maxima and minima have index +1 while saddle points have index -1.}

In any case, it is elementary to classify all possible combinations of non-degenerate critical points compatible with the formula \eqref{eq:index} \emph{and} with Bezout's theorem, which provides the maximal number of real critical points\footnote{Note that while the characteristic is an intrinsic property of the sphere we work on, the bound provided by Bezout's and by Rohrl's theorems depends on the degree of the mapping.}. {In view of the special nature of the potential considered here (homogeneous of degree $n$), it is either even or odd, depending on the parity of $n$, hence all critical points are dual to each other under reflection.} Thus in the odd case there will be as many maxima as minima while in the even one these numbers will necessarily be even, as well as (in all cases) the number of saddles of any given index.

For example, in \cite{GV2} it is argued that in the case of a cubic potential (and hence a quadratic gradient mapping) in three-dimensional space (and hence a two-dimensional ambient sphere) with non-degenerate extremals (but possibly degenerate saddles) only the possibilities listed in Table \ref{tab:TableI}  arise.\footnote{It should be noted that not all of these are realized when we consider, as the Physics of the problem studied in \cite{GV1,GV2} requires, \emph{completely traceless} tensors. See also the discussion in \cite{Wal} in this respect.}

Here we are more interested in even degree, and especially in \emph{quartic} potentials, and hence cubic gradient mappings. In this case, already for ambient space ${\bf R}^3$ we get up to $(3^3-1)  = 26$ critical points and a complete classification would make little sense. We remark, however, that Table  \ref{tab:TableI} is still valid in that it concerns topological features; on the other hand, in this context it does \emph{not} provide a complete classification of the possible situations, but only of those with no more than four maxima or minima.

We also note that when the ambient space is ${\bf R}^2$, and hence the relevant sphere is just a circle $S^1$, then we always have as many maxima as minima, whose number is of course limited by the degree of the potential; and of course no saddles. E.g. for $p=3$ and $q=2$ we have at most four maxima and four minima.

Some simple Examples will be considered in detail in Section \ref{sec:ExaI}.

\begin{table}
  \centering
  \begin{tabular}{|c|c|c|c|c||r|}
  \hline
  Max & Min & $S_1$ & $S_2$ & $S_3$ & NCP \\
  \hline
  1 & 1 & 0 & 0 & 0 & 2 \\
  2 & 2 & 2 & 0 & 0 & 6 \\
  3 & 3 & 4 & 0 & 0 & 10 \\
  3 & 3 & 0 & 2 & 0 & 8 \\
  4 & 4 & 6 & 0 & 0 & 14 \\
  4 & 4 & 2 & 2 & 0 & 12 \\
  4 & 4 & 0 & 0 & 2 & 10 \\
  \hline
\end{tabular}
  \caption{Different possibilities for the number and type of critical points in the case of a cubic potential in three dimensions; here ``Max''and ``Min'' represent the number of maxima and minima, while ``$S_k$''  represents the number of saddle points of index $- k$. Finally, ``NCP'' is the total number of critical points.}\label{tab:TableI}
\end{table}

\section{Example. Invariant lines; existence and number of higher order normal modes}
\label{sec:ExaI}

As the simplest possible example of the situation we have been
studying, we consider a point particle of mass $m=1$ in ${\bf
R}^2$ (with cartesian coordinates $x,y$) evolving under the action
of a quartic potential (a cubic one would not satisfy the requirement that the origin is a stable equilibrium); in order to reduce the complexity of the potential (and hence of the analysis) we assume it depends on $x$ and $y$ only through their squares. That is,  \beql{eq:V} V(x,y) \ = \ a \, x^4 \ + \ b \, y^4
\ + \ 2 \, c \, x^2\, y^2 \ . \eeq Note that this is symmetric
under $Z_2 \times Z_2$ (these acting as reflections in $x$ and in
$y$); if we wish to require this to be also invariant under the
$Z_2$ involution exchanging $x$ and $y$ then we should require
$b=a$. In this case it would be convenient, with a suitable
redefinition of constants, to rewrite this as \beql{eq:Vsym}
V(x,y) \ = \ \a \ (x^2 +y^2)^2 \ + \ \b \, x^2 \, y^2 \ . \eeq

We will refer to these cases as the \emph{lower symmetry} and the \emph{higher symmetry} cases respectively.

\subsection{The higher symmetry case}

We will first consider the case where the system is $Z_2 \times Z_2 \times Z_2$ symmetric, i.e. the potential is in the form \eqref{eq:Vsym}.

Note that in order to have a stable point at the origin, the
coefficients $\alpha$ should be positive; we will assume this to be the case from now on.

Moreover, we can always rescale $V$ (which amounts to a rescaling of time) and choose $\a =1$. Then in order to have a minimum at the origin we can ask $\b > - 4$; this is obtained by looking at the behavior of the potential $V_m (x) := V(x,mx)$ along all lines $y = m x$ (including the $m=\infty$ case, i.e. the $y$ axis). We will moreover assume $\b \not= 0$ to avoid the fully degenerate case with rotational symmetry.

When we pass to polar coordinates (we consider $\theta \in
(- \pi , \pi]$, and of course $\rho \in [0,\infty)$)
$$ x \ = \ \rho \ \cos (\theta) \ , \ \ \ y \ = \ \rho \ \sin (\theta ) $$ and
constrain the potential on the unit circle $\rho^2 =x^2+y^2 = 1$,
call it $W (\theta) $, we get
\beql{eq:W1} W \ = \frac18 \ \[ 8 \ - \ \b \, \cos (4 \, \theta ) \ + \ \b \] \ . \eeq
Thus we have
$$ \frac{d W}{d \theta} \ = \ \frac12 \ \b \ \sin (4 \, \theta ) \ ; $$
critical points are obtained for
$$ \theta \ = \ \theta_k \ := \ k \ \frac{\pi}{4} \ , \ \ \ |k| \le 4 \ . $$
That is, we get \emph{eight} critical points on the unit circle, corresponding to \emph{four} invariant lines; this applies for \emph{any} nonzero value of $\b$ (as already remarked $\b=0$ is the fully rotationally invariant and hence infinitely degenerate case; we excluded this from our considerations).

\begin{figure}
  \includegraphics[width=200pt]{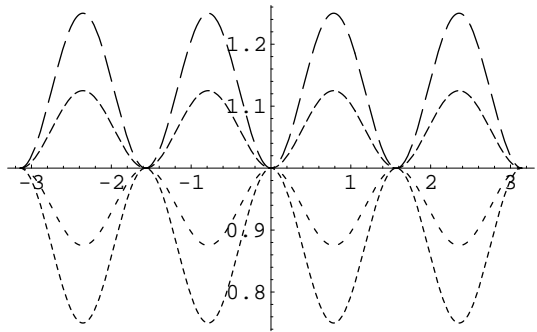}\\
  \caption{The potential $W (\theta)$ as in eq.\eqref{eq:W1} for different values of $\b$; here $\b =-1, -0.5,0.5,1$. The exchange of stability takes place at $\b=0$.}\label{fig:plotsymm}
\end{figure}

The stability of the critical points is controlled by
$$ \[ \frac{d^2 W}{d \theta^2} \]_{\theta_k} \ = \ 2 \ \b \ \cos (4 \, \theta_k ) \ . $$
Thus we have a bifurcation at $\b = 0$; for this value of $\b$ all the stabilities are exchanged. In particular, for $\b < 0$ the lines identified by $\theta = \pm \pi/4$ are stable and the axes ($\theta = 0,\pi/2$) are unstable, while for $\b>0$ the lines identified by $\theta = \pm \pi/4$ are unstable and the axes are stable; see Figure \ref{fig:plotsymm}.\footnote{We stress that here the stability means stability for the potential restricted to the unit circle, while the origin is \emph{always} stable for the full two-dimensional potential $V(x,y)$.}

The effective potentials along the invariant lines, call these as $V_\theta$ where $\theta$ is the (invariant) angular coordinate in the $(x,y)$ plane, is always of the form $V_\theta (r) = c_\theta \, r^4$ with $c_\theta$ a constant.
We actually get (recall we assumed $\b > - 4$)
$$ c_0 \ = \ c_{\pi/2} \ = \ 1  \ ; \ \ c_{\pi/4} \ = \ c_{- \pi/4} \ = \ 1 \ + \ \b/4 \ . $$ This is coherent, of course, with the stability of the equilibrium point at the origin.

We stress that, apart from the stability exchange for $W (\theta)$, in this case there are no qualitative changes as the parameters (which in this case means just the parameter $\b$) are varied: we always have four critical lines and hence four normal modes; two of them are stable and two of them unstable.

\begin{figure}
\centering
  \begin{tabular}{cc}
  \includegraphics[width=100pt]{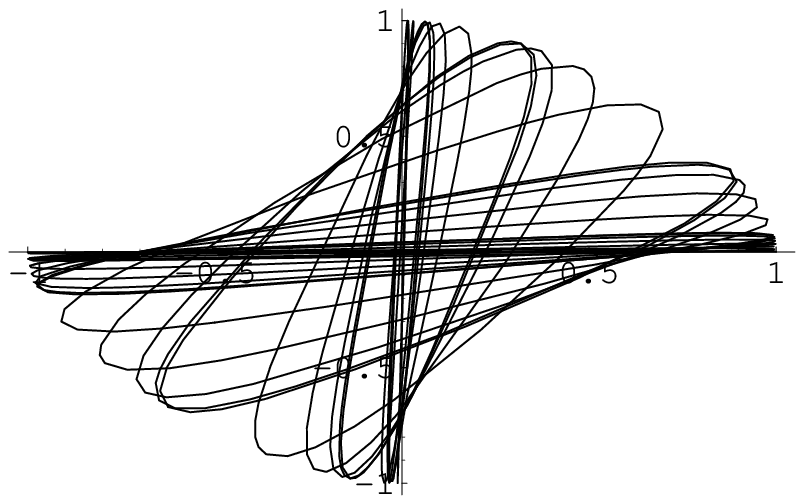} & \includegraphics[width=100pt]{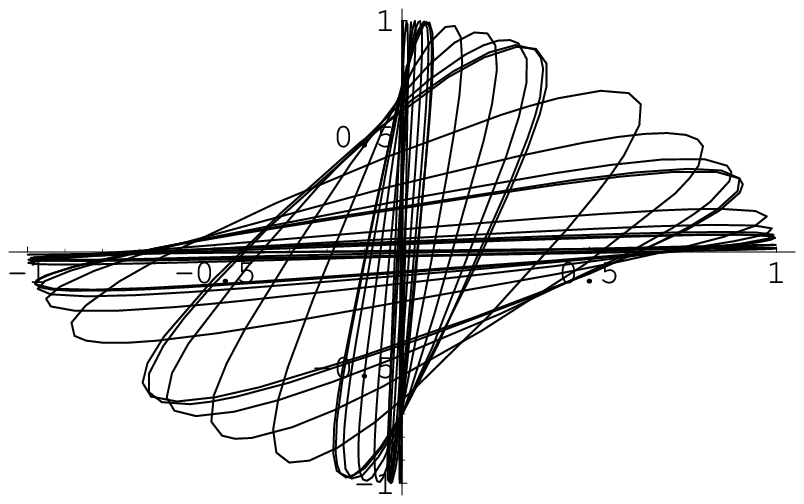}\\
  $(a)$ & $(b)$ \\
  \includegraphics[width=100pt]{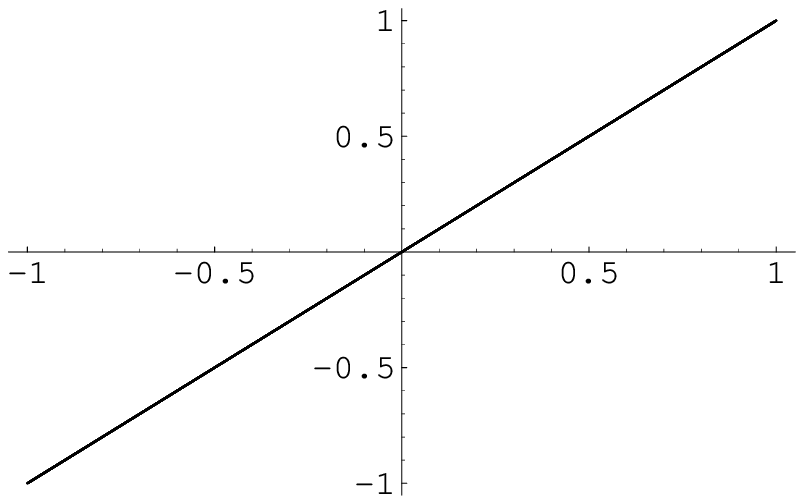} & \includegraphics[width=100pt]{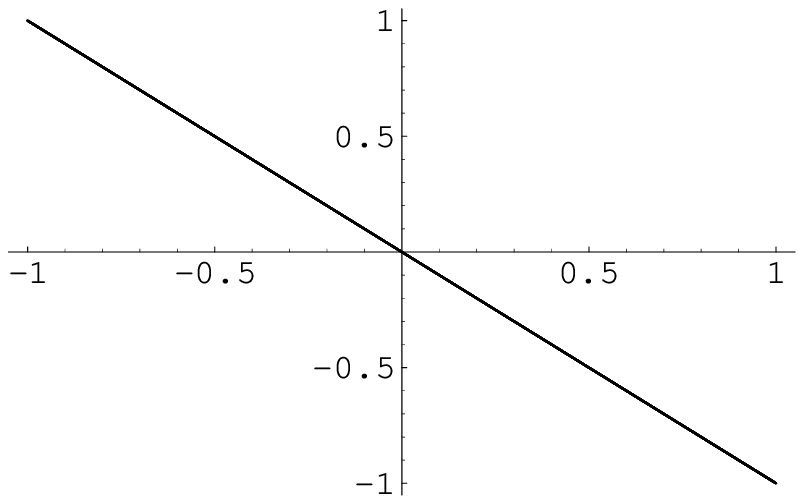}\\
  $(c)$ & $(d)$ \end{tabular}
  \caption{Numerical integration of the motion generated by the potential \eqref{eq:Vsym} with the choice $\b = 1$ for initial conditions near to normal modes. In all cases, initial data correspond to zero speed and position at $r=1$ along eigenvectors, with an offset of 0.001 from the latter. The simulation show the outcome, for $t \in (0,100)$, for initial data: (a) near the eigenvector $\theta =0$, (b) near the eigenvector $\theta = \pi$, (c) near the eigenvector $\theta = \pi/4$, (d) near the eigenvector $\theta = - \pi/4$.}\label{fig:ExaIsymm}
\end{figure}

Some numerical simulations of this dynamics for initial data near to the (higher order) normal modes are shown in Figure \ref{fig:ExaIsymm}; they confirm stability as discussed above.\footnote{Here we show the situation for $\b>0$; simulations for $\b< 0$ would also confirm our discussion, and are not shown for the sake of brevity.}

\subsection{The lower symmetry case}

We will now consider the general form (which has only $Z_2 \times Z_2$ as symmetry), i.e. the potential \eqref{eq:V}. With no loss of generality, we can assume $a > b$ (if not, just switch $x$ and $y$). Here again stability requires that both $a$ and $b$ are positive; by a rescaling we can set $a=1$, and hence $0 < b < 1$, and deal with the potential
\beql{eq:Vs} V(x,y) \ = \ x^4 \ + \ \a \, y^4
\ + \ 2 \, \b \, x^2\, y^2 \ , \eeq where $0 < \a < 1$.

By looking at this along the line $y =m x$ we get
$$ V_m (x) \ := \ V(x,mx) \ = \ (1 + 2 \b m^2 + \a m^4 ) \ x^4 \ ; $$
stability of the origin requires therefore that $$ 1 \ + \ 2 \, \b \, m^2 \ + \ \a \, m^4 \ > \ 0 $$ for all choices of $m$, and this implies
$$ \b \ > \ - \ \sqrt{\a} \ ,  $$ which we assume from now on.

After passing to polar coordinates, the restriction of $V$ given by \eqref{eq:Vs}  to the unit sphere reads
\beql{eq:W} W (\theta) \ = \ \cos^4 (\theta) \ + \ \a \ \sin^4 (\theta) \ + \ 2\, \b \, \sin^2(\theta) \cos^2 (\theta) \ , \eeq and from this we get at once
$$ \frac{d W}{d \theta} \ = \ - \ \[ 1 \, - \, \a \, + \, (1 + \a - 2 \b) \, \cos (2 \theta) \] \ \sin (2 \theta) \ . $$
Thus critical points of $W$ are identified either by $\sin (2 \theta) =0$, i.e. by $\theta = 0 , \pm \pi/2 , \pm \pi$; or by
$$ \cos (2 \theta) \ = \ \frac{\a -1}{\a+1 - 2 \b} \ . $$
Solutions to this equation exist only {in the regions} $\b < \a$ and $\b > 1$ (recall that we assumed $0< \a <1$), but not for $\a < \b < 1$.

Thus we conclude that $\b$ lies in the range $\b > -\sqrt{\a}$, and that for $\b$ taking the values $\b = \a$ and $\b = 1$  there are bifurcations changing the number of critical points for $W$, i.e. of invariant lines for our potential $V(x,y)$.

More precisely, we always have four critical points at $\theta = 0 , \pm \pi/2 , \pi$; moreover for suitable $\b$ we also have four more critical points at
$$ \theta \ = \ \pm \ \frac12 \ \arccos \[ \frac{1 - \a}{1 + \a - 2 \b} \] \ := \ \pm \theta_* \ . $$

\begin{figure}
\begin{tabular}{ccc}
  \includegraphics[width=75pt]{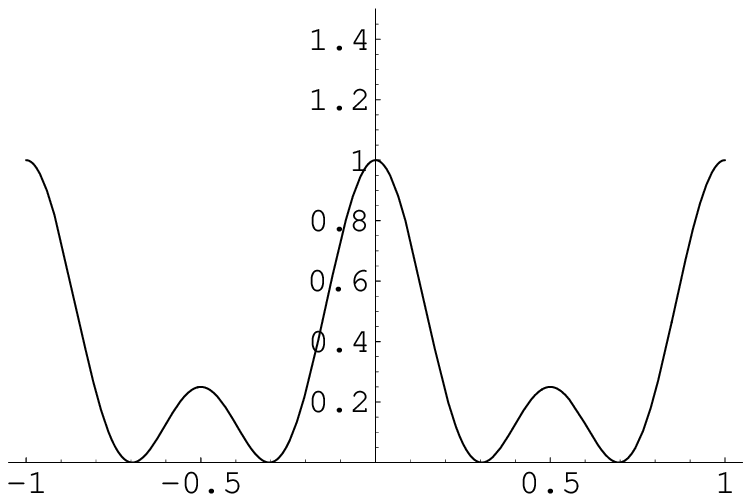} &
  \includegraphics[width=75pt]{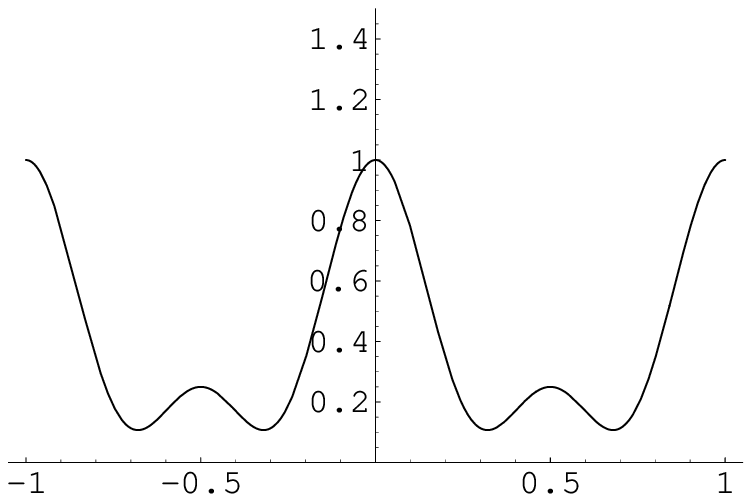} &
  \includegraphics[width=75pt]{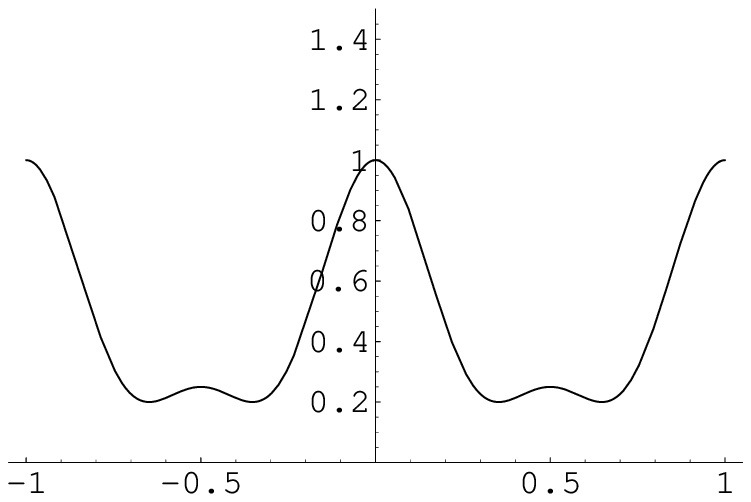} \\
  $\b = - 1/2$ & $\b = -1/4$ & $\b = 0$ \\
   & & \\
  \includegraphics[width=75pt]{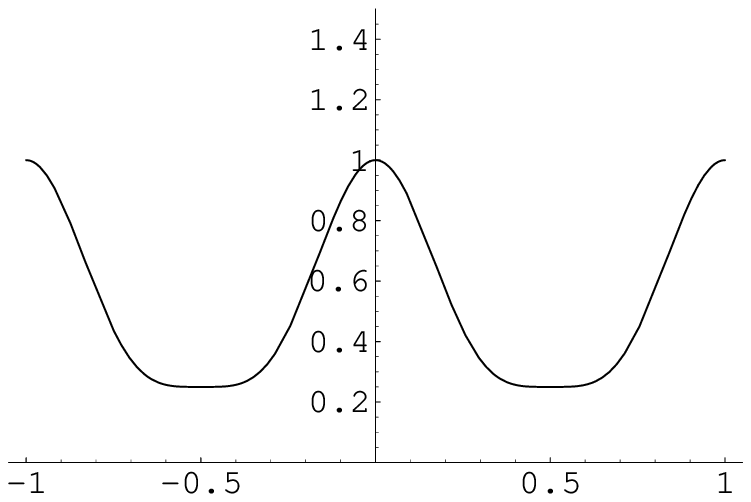} &
  \includegraphics[width=75pt]{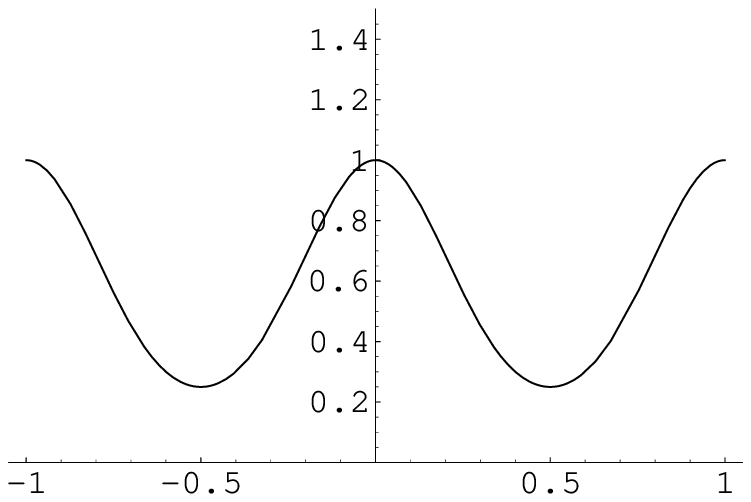} &
  \includegraphics[width=75pt]{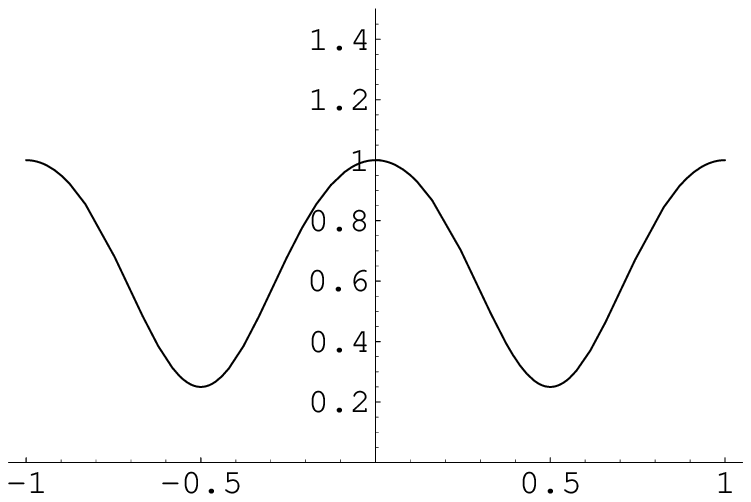} \\
  $\b = 1/4$ & $\b = 1/2$ & $\b = 3/4$ \\
  \includegraphics[width=75pt]{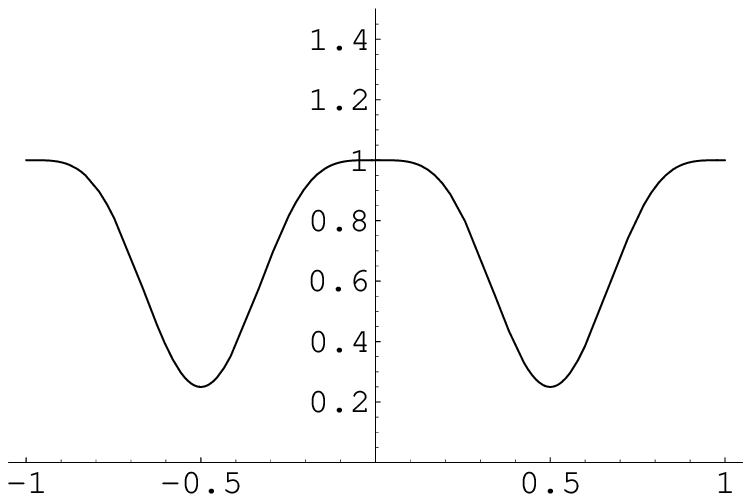} &
  \includegraphics[width=75pt]{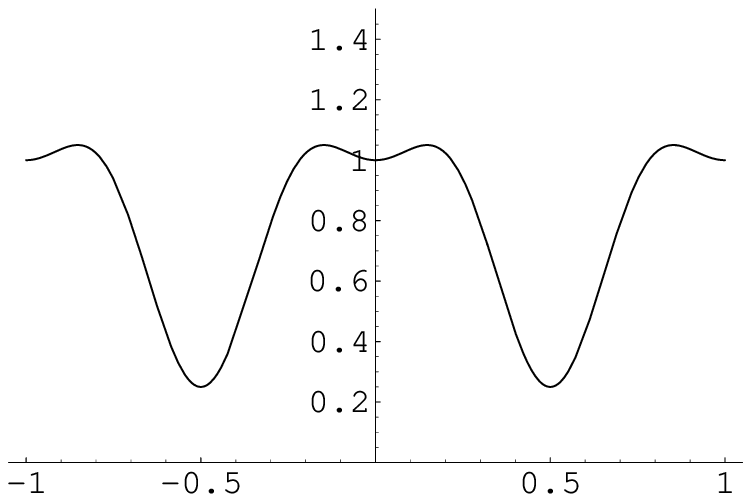} &
  \includegraphics[width=75pt]{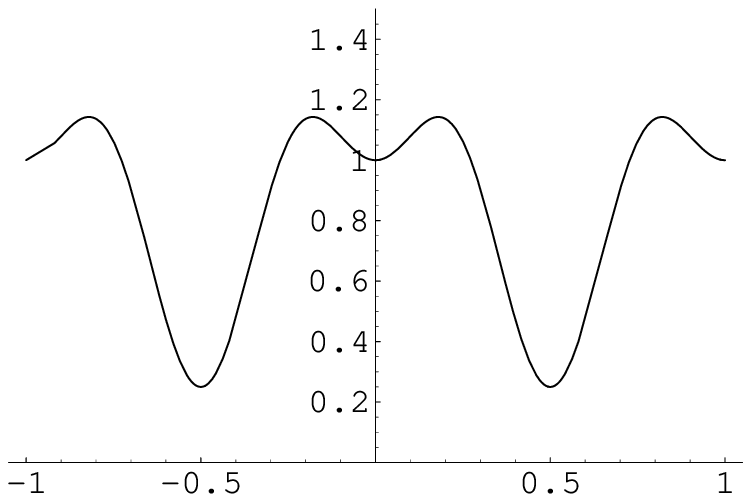} \\
  $\b = 1$ & $\b = 5/4$ & $\b = 3/2$ \\
 \end{tabular}
  \caption{The potential $W(\theta)$, see \eqref{eq:W}, for $\a = 1/4$ and various choices of $\b$. Here $\theta$ is measured in units of $\pi$.}\label{fig:plotNS}
\end{figure}

That is, we have either eight or four critical points for the potential $W$ restricted on the sphere, corresponding to \emph{four} or \emph{two} critical lines and hence normal modes, depending on the value of $\b$. See Figure \ref{fig:plotNS}.

The stability of critical points is controlled by
$$ \frac{d^2 W}{d \theta^2} \ = \ - \, 2 \ \[ (1 - \a) \ \cos (2 \, \theta) \ + \ (1 + \a - 2 \b) \ \cos (4 \theta ) \] \ ; $$
in particular at the various critical points identified above we have
\begin{eqnarray*}
(d^2 / d \theta^2)_{\theta = 0} &=& (d^2 / d \theta^2)_{\theta = \pi} \ = \ - \ 4 \ (1 - \b) \ ; \\
(d^2 / d \theta^2)_{\theta = \pi/2} &=& (d^2 / d \theta^2)_{\theta = -\pi/2} \ = \ - \ 4 \ (\a  - \b) \ ; \\
(d^2 / d \theta^2)_{\theta = \theta_*} &=& (d^2 / d \theta^2)_{\theta = - \theta_*} \ = \ 8 \ \[ \frac{(\a -\b) \ (1 - \b)}{1 + \a -2 \b} \] \ . \end{eqnarray*}

This immediately shows that the critical line corresponding to $\theta =0$ (equivalently, to $\theta = \pi$) undergoes a change of stability at $\b = 1$; and the critical line corresponding to $\theta = \pi/4$ (equivalently, to $\theta = 3 \pi /4$) undergoes a change of stability at $\b = \a$. As for the critical lines corresponding to $\theta = \pm \theta_*$, existing in the regions $\b < \a$ and $\b > 1$, noting that $\a < (1+\a)/2 < 1$, we have that these are stable for $\b < \a$ and unstable for $\b >1$. (Numerical simulations, not shown for the sake of brevity, confirm again our analysis.)

%\newpage

\section{Discussion and Conclusions}

In this note we have considered natural Hamiltonians for point particle, $H = K+V$, for which the dynamics in local coordinates is simply $\ddot{x} = - \nabla V$. We have considered the neighborhood of a stable equilibrium $x_0$, and studied the case where the Taylor expansion of $V$ starts with terms of order $k$ higher than two (we have of course given special attention to the case $k=4$).

We have discussed how the notion of \emph{normal modes} is modified in this case; this is based essentially on known results concerning \emph{eigenvalues of (homogeneous) tensors}. It results that higher order normal modes do exist, but while some of their properties extend from the standard (i.e. $k=2$) case to the present one, other do not. In particular, their number is not fixed and can exceed the dimensionality of the ambient space; moreover the most general dynamics near the equilibrium is in general\footnote{It is conceivable that one can build systems for which a \emph{nonlinear superposition principle} \cite{NLSP1, NLSP2} holds. We did not discuss this point here.} \emph{not} a superposition of normal modes, at difference with the standard case.

It should be mentioned that a question which arise naturally has not been studied here, and should be considered in the future. This is of course the \emph{persistence} of these higher order normal modes under perturbations, i.e. when one considers also higher order terms in the series expansion for the potential around the equilibrium point. We recall that for standard normal modes a theory of persistence exists \cite{NLNM1,NLNM2}; this is based on variational analysis and guarantees persistence of some of the normal nodes under certain conditions. Apart from an extension to the new higher order normal modes along this line of attack, one could also consider an approach based on the theory of Poincar\'e-Birkhoff \emph{normal forms} \cite{ArnGM,Elp} (or some generalization thereof).

It is worth noting that for the examples discussed in Section \ref{sec:ExaI}, persistence of some modes can be guaranteed by symmetry arguments if the symmetry is preserved by higher order perturbations: In both cases, the potential then admits symmetries sending $x\mapsto x,\,y\mapsto -y$, resp.  $x\mapsto -x,\,y\mapsto y$ (analogously for the time derivatives), hence the fixed point spaces of these symmetries are necessarily invariant. The fixed point spaces correspond to the cases $\theta= 0,\,\pm\pi$ above. In the higher symmetry case one also has symmetries exchanging $x$ and $y$, resp. $x$ and $-y$, with invariant fixed point spaces corresponding to $\theta=\pi/4$ resp.\ $\theta=3\pi/4$.

We plan to tackle the persistence problem generally in a later publication.

\end{document}